\documentclass[12pt]{article}
\usepackage{amsmath, amssymb, amsthm, amsfonts}
\usepackage{graphicx}
\usepackage{xcolor}
\usepackage{enumitem}
\usepackage{hyperref}
\usepackage{tikz}
\usepackage{tikz-cd}
\usepackage{tikz-cd}
\usepackage{pgfplots}
\pgfplotsset{compat=1.18}
\usetikzlibrary{arrows, positioning, calc, patterns, arrows.meta}

\newtheorem{theorem}{Theorem}[section]
\newtheorem*{theoremA}{Theorem A}
\newtheorem*{theoremB}{Theorem B}
\newtheorem{lemma}[theorem]{Lemma}
\newtheorem{proposition}[theorem]{Proposition}
\newtheorem{corollary}[theorem]{Corollary}
\theoremstyle{definition}
\newtheorem{definition}[theorem]{Definition}
\theoremstyle{remark}
\newtheorem{remark}[theorem]{Remark}

\title{Geometry of Hypersurfaces with Isolated Singularities}
\author{Jiayi Hu, Fengyang Wang, Xinlang Zhu}
\date{}

\begin{document}

\maketitle
\begin{abstract}
This paper explores the Fano variety of lines in hypersurfaces, particularly focusing on those with mild singularities. Our first result explores the irreducibility of the variety $\Sigma$ of lines passing through a singular point $y$ on a hypersurface $Y \subset \mathbb{P}^n$. Our second result studies the Fano variety of lines of cubic hypersurfaces with more than one singular point, motivated by Voisin's construction of a dominant rational self map. 
\end{abstract}


\section{Introduction}
\label{sec:introduction}

The study of the Fano variety of lines in a projective hypersurface plays an important role in algebraic geometry. A classical result of Cayley and Salmons states that a smooth cubic surface contains precisely $27$ lines. The Fano variety of lines in a cubic threefold is thoroughly studied in~\cite{ClemensGriffiths}. Moreover, the Fano variety of lines in a cubic fourfold is established as a hyper-Kähler fourfold in~\cite{BeauvilleDonagi}. 

A prevalent method for studying smooth hypersurfaces, or their associated moduli spaces, involves their degeneration to cases with singularities. Noether-Lefschetz theory serves as an important example~\cite[Chapter 2]{VoisinHodgeII}. Additionally, intersecting a smooth hypersurface with a linear subspace can naturally induce singularities~\cite{Bai, ClemensGriffiths}. This motivates the investigation of the Fano variety of lines in a hypersurface with mild singularities, leading to our initial result summarized below. 

\begin{theoremA}
    Let $Y\subset \mathbb P^n$ be a hypersurface of degree $d$ admitting a point $y\in Y$ with multiplicity $m$ as the only singularities. Assume that $(Y, y)$ is general in the moduli space of all such pairs. Let $\Sigma$ be the variety of lines in $Y$ passing through $y\in Y$. Then\\
    (i) If $d\leq m + n - 2$, then $\Sigma$ is a smooth variety of dimension $m + n - 2 - d$.\\
    (ii) If $d = m + n -2$, then $\Sigma$ has $\frac{d!}{(m - 1)!}$ points.\\
    (iii) If $d < m + n - 2$, then $\Sigma$ is irreducible.
\end{theoremA}

Consider Theorem A in the following specific scenarios: 

\begin{itemize} 
\item For $d = 3, n = 3, m = 2$, Theorem A shows that $\Sigma$ is a finite set consisting of $3!/1! = 6$ points. This corresponds to the classical understanding that a cubic surface with a single simple node as a singularity has exactly $6$ lines passing through this point. 
\item For $d = 3, n = 4, m = 1$, Theorem A shows that $\Sigma$ is a finite set consisting of $3! / 0! = 6$ points. This corresponds to the classical result~\cite{ClemensGriffiths} that there are exactly $6$ lines passing through a general point on a smooth cubic threefold. 
\end{itemize}

A direct generalization of the Fano variety of lines in a hypersurface is the scheme of linear subspaces in a complete intersection. Some basic yet important properties are extensively studied in~\cite{DM}. Motivated by the fundamental example of Beauville-Donagi~\cite{BeauvilleDonagi}, the following Calabi-Yau manifolds are constructed in~\cite{KCorr}, as important examples to illustrate intriguing conjectures in algebraic geometry such as the Kobayashi hyperbolicity conjecture and the Voisin's conjecture on K-trivial varieties. Let $r, n$ be positive integers satisfying the numerical condition
\[
 n + 1 = \binom{ r + 3 }{2}.
\]
Let $Y\subset \mathbb P^n$ be a general smooth cubic hypersurface of dimension $n - 1$. Consider the scheme $X$ of linear subspaces of dimension $r$ inside $Y$. Then it is proven~\cite{KCorr} that $X$ is a K-trivial variety endowed with a dominant rational map $\Psi: X\hookrightarrow X$, now called the Voisin map, defined based on the following important construction. Let $x\in X$ be a general point corresponding to an $r$-linear subspace $P_x\subset Y$. There exists a unique linear subspace $H_{x, r + 1}\subset \mathbb P^n$ of dimension $r + 1$ that is tangent to $Y$ along $P_x$. The intersection of $H_{x, r+1}\cap Y$ should be a cubic hypersurface in $H_{x, r+1}$, but since the intersection already contains a double $P_x$, there remains only a residual linear subspace $P_{x'}$ in the intersection. The Voisin map is defined as sending $x$ onto the point $x'$ corresponding to the linear subspace $P_{x'}$. The geometry properties of the Voisin map have been extensively studied in~\cite{Bai}, especially in the special case $r=2$. Motivated by the construction in~\cite{Bai}, we prove the following result.

\begin{theoremB}
    With notations as above, let $H_{2r + 1}\subset \mathbb P^n$ be a general linear subspace containing $H_{x, r + 1}$. Let $Y_{2r} = H_{2r + 1} \cap Y$. Then\\
    (i) $Y_{2r}$ has $2^r$ simple double points as the only singularities.\\
    (ii) Let $y\in Y_{2r}$ be one of the singular point. Let $\Sigma_y$ be the variety of lines in $Y_{2r}$ passing through $y$. If $r\geq 2$, then $\Sigma_y$ is irreducible of dimension $2r - 2$.
\end{theoremB}

Theorem B generalizes~\cite[Lemma 3.36, Lemma 3.37]{Bai}.

\section{Preliminaries}
\label{sec:preliminaries}

In this section, we will go through some basic concepts and results that are fundamental to this work.
\subsection{Definitions and Notation}

Let $X$ be a complex algebraic variety. The relation between the line bundles and divisors serves as an important tool to study the geometry of $X$. Let us recall these classical notions here for the readers' convenience.
\begin{definition}[Line Bundles]
    A line bundle on $X$ is a morphism $\pi : L \to X$, such that there exists an open cover $\{U_\alpha\}_{\alpha \in \Lambda}$ of $X$, such that on each $U_{\alpha}\subset X$,

    \begin{center}
    \begin{tikzcd}
        \pi:\quad L_{U_{\alpha}} := \pi^{-1}(U_\alpha) \arrow[r, "\pi"] \arrow[d, "\cong", "\varphi_\alpha"'] & U_{\alpha} \\
        \mathbb{C}\times U_{\alpha} \arrow[ur]
    \end{tikzcd}    
    \end{center}
    where $\varphi_{\alpha}$ may not be uniquely determined by $L$. In fact, $\varphi_{\alpha}$ is determined up to the multiplication of an invertible function on $U_{\alpha}$.
    Furthermore, for any $\alpha,\beta \in \Lambda$, the following must be satisfied
    \begin{center}
    \begin{tikzcd}
        (U_{\alpha} \cap U_{\beta}) \times \mathbb{C} & \arrow[l, "\varphi_{\beta}"', "\cong"]
        L_{U_{\alpha} \cap U_{\beta}} \arrow[r, "\varphi_{\alpha}", "\cong"'] &
        (U_{\alpha} \cap U_{\beta}) \times \mathbb{C}.
    \end{tikzcd}
    \end{center}
Finally, we require that the map $g_{\alpha\beta}: U_\alpha \cap U_\beta \to \mathbb C^*$ induced by the composition of the above maps is a morphism.
\end{definition}

Let $s$ be a non-zero rational section of $L$, then the zero locus of $s$ is a finite union of irreducible subvarieties $\bigcup_i Z_i$ of codimension $1$ of $X$. Similarly, the poles of $s$ is also a finite union of irreducible subvarieties $\bigcup_j Z_j'$. Let $m_i$ (resp. $m_j'$) be the multiplicity of the section $s$ on the zeros $Z_i$ (resp. poles $Z_j'$). The formal sum $\sum_i m_iZ_i - \sum_j m_j'Z_j'$ is called the divisor associated to the section $s$, written as $\mathrm{div}(s)$. 

In general, a divisor $D$ of $X$ is a formal sum of the following form
\[
D = \sum_i n_iZ_i,
\]
where $n_i\in\mathbb Z$ and $Z_i$ is an irreducible subvariety of codimension $1$ in $X$. One of the most important equivalence relations of the divisors is the linear equivalence.

\begin{definition}[Linear Equivalence]
    Two divisors $D_1, D_2$ are called linearly equivalent, if $D_1 - D_2$ is a principal division, i.e. there exists a rational function $f$ such that 
    \[
    D_1-D_2=\operatorname{div}(f)
    \]
    where $\operatorname{div}(f)$ is the divisor associated to the rational function $f$ (viewed as a rational section of the trivial line bundle).
\end{definition}

A divisor $D = \sum_i n_i Z_i$ ($n_i \neq 0$ for all $i$) is called \emph{effective} if all the coefficients $n_i$ are nonnegative. The \emph{support} of $D$ is $\mathrm{supp}(D) := \bigcup_i Z_i$ the union of all its irreducible components. 

\begin{definition}[Linear Systems and Base Loci]
   Let $D$ be a divisor of $X$. The linear system associated with $D$ is 
    \[
        |D|:=\{ D' \text{ effective divisor } |D'
        \textrm{ is linearly equivalent to } D \}.
    \]
    A base point $x\in X$ of $|D|$ is a point that is contained in any member of the linear system $|D|$.
    The base locus of $|D|$ is the set of base points of $|D|$. Equivalently, the base locus of $|D|$ is the intersection of the supports of all members of $|D|$.
\end{definition}

\begin{definition}[Multiple Point $\mathbb{A}^n$]
    Let $Y\subset \mathbb A^n$ be a hypersurface defined by a function $f$. A point $y\in Y$ is called a multiple point of multiplicity $m$, if in the Taylor expansion of $f$ at $y$ has no terms of degree lower than $m$, but has degree $m$ term.
\end{definition}
\begin{definition}[Multiplicity in $\mathbb{P}^n$]
    Let $Y\subset \mathbb P^n$ be a hypersurface. The multiplicity of $a\in Y\subset \mathbb{P}^n $ is defined to be the multiplicity of $a$ viewed in some affine chart containing $a$. 
\end{definition}
\begin{definition}[Degree of Subvariety]
    Let $Z\subset \mathbb{P}^N$ be a subvariety of dimension $n$. The degree of $Z$ is defined to be the number of intersection point of $Z$ with a general linear subspace of codimension $n$.
\end{definition}


\subsection{Classical Results}
In this section, we review two key classical results used in the proofs of Theorem A and Theorem B.
\begin{theorem}[Bertini theorem]\label{bertini1} Let $X$ be an algebraic variety over $\mathbb{C}$. Let $|D|$ be a linear system of $X$. Then for a general element $Y\in |D|$, $Y$ is smooth outside the base locus of $|D|$ and the singular locus of $X$.
\end{theorem}
A special case of the Bertini theorem is the following
    \begin{corollary}\label{bertini2}
     Let $Y\subset \mathbb{P}^N$ be a smooth subvariety. Let $f\in \mathbb{C}[x_0,...,x_N]$ be a homogeneous polynomial such that $f|_Y$ is not constantly zero. Then $V(f)\cap Y$ is of dimension $N-1$. Furthermore, if $f$ is general, then $V(f)\cap Y$ is smooth.   
\end{corollary}
The linear system $|D|$ is a projective space. A element $f\in |D|$ is called general if the element $f$ is taken in a nonempty Zariski open subset of $|D|$. A proof of the Bertini theorem can be found in~\cite{Bertini}.

Another important result is the Fulton-Hansen connectedness theorem.
\begin{theorem}[Fulton-Hansen connectedness theorem]\label{connectedness}
    Let $V,W\subset \mathbb{P}^n$ be irreducible subvarieties such that 
    \[
    \operatorname{dim}V+\operatorname{dim}W>n.
    \]
    Then $V\cap W$ is connected.
    \end{theorem}
One can refer to~\cite{FultonHansen} for a proof of this result.
 \begin{remark} Note that when $\dim V + \dim W = n$, the intersection $V \cap W$ maybe not connected. For a counter-example, we consider the circle $V$ defined by \( V = \{(x : y : z) \in \mathbb{P}^2 \mid x^2 + y^2 = z^2\} \) $\mathbb{P}^2$ and a hyperbola $W$ defined by \( W = \{(x : y : z) \in \mathbb{P}^2 \mid xy = z^2\} \). In this case, we have \(\dim V + \dim W = 2\). But the intersection $V\cap W$ is a set of isolated points. In particular, $V\cap W$ is not connected.
\end{remark}


\section{The Geometry of Lines on Hypersurfaces}
\label{sec:geometry-of-lines}

We study the geometry of lines on hypersurfaces in this section and prove Theorem A. The following classical result is used in our argument and we write a proof of it for completeness.
\begin{proposition}\label{PropClassicalResult}
    Let $f_1, \ldots, f_k\in \mathbb C[x_0, \ldots, x_N]$ be general homogeneous polynomials of degree $d_1, \ldots, d_k$, respectively. Let $Z\subset \mathbb P^N$ be the zero locus of $f_1, \ldots, f_k$. If $k\leq N$. Then $Z\subset \mathbb P^N$ is a smooth subvariety of dimension $N - k$ and of degree $d_1d_2\ldots d_k$. If furthermore, $k < N$, then $Z$ is irreducible.
\end{proposition}


\begin{proof}
  
    It suffices to show that requiring $f_k|_Y$ to be not constantly $0$ is a general condition. Granting it, by the theorem of Bertini (Corollary \ref{bertini2}), for any general $f_1,\dots,f_{k-1}\in\mathbb{C}[x_0,...,x_N]$, the variety
    \[Y:=V(f_1)\cap \dots \cap V(f_{k-1})\subset \mathbb{P}^N\]
    is smooth of dimension $N-k+1$. Then for any general homogeneous polynomial $f_k\in \mathbb{C}[x_0,\dots,x_N]$ and $f_k|_Y$ not constantly equal to $0$, the variety $Y\cap V(f_k)$ is smooth of dimension $N - k$.
    
    Let us now prove that requiring $f_k|_Y$ to be not constantly zero is a general condition. We define
    \begin{align*}
    \operatorname{ev}_Y:\mathbb{C}_{d_{k}}[x_0,...,x_N] & \longrightarrow \{\text{Rational functions on $Y$ of degree $d_k$}\}\\
        f &\longmapsto f|_Y
    \end{align*}
    The map $\operatorname{ev}_Y$ is a linear map and is continuous under the Zariski topology. Hence, the set 
    \[\operatorname{ev}_Y^{-1}(\{0\})^C= \left \{f_k\in \mathbb{C}_{d_k}[x_0,...,x_N]\ \mid\ f_k|_Y \text{ not constantly zero}\right\}\]
 is Zariski open in $\mathbb{C}_{d_k}[x_0,...,x_N]$. Clearly, $\operatorname{ev}_Y^{-1}(\{0\}^C\neq \emptyset$. In fact, for any $y=[y_0:\cdot\cdot\cdot:y_N]\in Y \subset \mathbb{P}^N$ (may assume $y_0\neq 0$ without loss of generality), we can find a homogeneous polynomial $f_k\in \mathbb{C}_{d_{k}}[x_0,...,x_N]$ of degree $d_k$, such that. $f_k(y)\neq 0$. Specially, we take
    \[f_k(x_0,...,x_N)=(x_1-\frac{y_1}{y_0}x_0)^{d_k}+x_0^{d_k}.\]
    From definition, the homogeneous polynomial $f_k|_Y$ is not constantly zero. If $k<N$, we have 
    $$\dim Y+\dim V(f_k)=N-k+1 + N - 1=2N-k>N.$$
 By \autoref{connectedness}, $Z=Y\cap V(f_k)$ is connected, and hence, $Z$ is irreducible since any smooth connected variety is irreducible.
\end{proof}

In the rest of this section, let $Y\subset \mathbb P^n$ be a general hypersurface defined by $f(x_0, \ldots, x_n)$ of degree $d$ admitting a point $y\in Y$ with multiplicity $m$ as the only singular point. We may assume that the coordinate of the singular point $y$ is $[1: 0: \ldots: 0]$ in $\mathbb P^n$. By the definition of the multiplicity of a singular point, we have 
\begin{lemma} The polynomial $f$ is of the following form:
    \[
    f(x_0, x_1, \ldots, x_n) = \sum_{i = m}^d x_0^{d - i} f_i(x_1, \ldots, x_n),
    \]
    where $f_i\in \mathbb C[x_1, \ldots, x_n]$ is a homogeneous polynomial of degree $i$.
\end{lemma}
\begin{proof}
    On the open subset $U_0:=\{ (1: x_1: \ldots: x_n) \in \mathbb P^n \} \cong \mathbb A^n$, the subvariety $Y\cap U_0$ is defined by a polynomial $g(x_1, \ldots, x_n) : = f(1, x_1, \ldots, x_n)$. Let $f_i(x_1, \ldots, x_n)$ be the degree $i$ homogeneous part of $g$. Then the Taylor expansion of $g$ at the point $y = (0, \ldots, 0)\in \mathbb A^n \cong U_0 \subset \mathbb P^n$ is nothing but $\sum_{i = 0}^d f_i(x_1, \ldots, x_n)$. The multiplicity of $y\in Y$ is $m$, so we have $f_i(x_1, \ldots, x_n) = 0 $ for any $i < m$ and $f_m(x_1, \ldots, x_n) \neq 0$. Hence, $g(x_1, \ldots, x_n) = \sum_{i = m}^n f_i(x_1, \ldots, x_n)$. Finally, 
    \[
    f(x_0, x_1, \ldots, x_n) = x_0^d\cdot g\left(\frac{x_1}{x_0}, \ldots, \frac{x_1}{x_0}\right) = \sum_{i = m}^d x_0^{d - i } f_i(x_1, \ldots, x_n),
    \]
    as desired.
\end{proof}

Let $H\cong \mathbb P^{n - 1}$ be the hyperplane in $\mathbb P^n$ defined by $x_0 = 0$. For the sake of convenience, we take the following identification.

\begin{lemma}\label{LmmLinesInProjectiveSpace}
    The set of lines in $\mathbb P^n$ passing through the point $y$ can be identified with $H$ by sending a line to the intersection point of the line and $H$.    
\end{lemma}

\begin{proof} If $\ell_{x}$ is a line in $\mathbb{P}^n$ passing through $y$ and another point $x=[x_0:...:x_n]$, then the intersection of $\ell_x$ with $H$ is the point $[0;x_1y_0-x_0y_1:...:x_ny_0-x_0y_n]$. The map from set of lines in $\mathbb P^n$ passing through the point $y$ to $H$ given by $\ell_{x}\mapsto[0:x_1y_0-x_0y_1:...:x_ny_0-x_0y_n]$ is clearly a bijection, hence the lemma follows.
\end{proof}
Let $L$ be a line passing through $y$ in $\mathbb P^n$ and let $y' = (0: x_1: \ldots: x_n)$ be the intersection point of $L$ and $H$. It is not hard to see that
\begin{lemma}\label{LmmExpressionDefiningSigmaInH}
    The whole line $L$ lies in $Y$ if and only if for any $i = m, m + 1, \ldots, d$, we have
    \[
    f_i(x_1, \ldots, x_n) = 0.
    \]
\end{lemma}
\begin{proof}
    The line $L$ passes through $y=(1:0:\dots:0)$ and $y'=(0:x_1:\dots:x_n)$, so $L$ is of the form 
    \[
    L=\{(u:vx_1:\dots:vx_n)\:|\:u,v\in \mathbb{C}\}.
    \]
    
    Assume that $f_i(x_1,...,x_n)=0$ for any $i=m, m+1,...,d$, from \mbox{Lemma 3.2.} we may write
    \begin{align*}
    f(x_0,x_1,\dots,x_n)=&x_0^{d-m}f_m(x_1,\dots,x_n)+\\
    &x_0^{d-m-1}f_{m+1}(x_1,\dots,x_n)+\dots+\\
    &x_0^0f_d(x_1,\dots,x_n).
    \end{align*}
    Then 
    \[
    f(x)=f(u,vx_1,\dots,vx_n)=\sum_{i = m}^d u^{d - i} v^i f_i(x_1, \ldots, x_n) = 0\]
    for all $x\in L$.
    This gives $L\subset Y.$

  Conversely, given that $L\subset Y$, then for any $x\in L$, we have
    \[
    0=f(x)=f(u,vx_1,\dots,vx_n)=u^{d-m}v^mf_m(x_1,\dots,x_n)+\dots+u^0v^df_d(x_1,\dots,x_n).
    \]
   Since $u,v$ are arbitrary, $f_i(x_1,\dots,x_n)=0$ for all $i=m,m+1,\dots,d$.    
\end{proof}

\begin{proof}[Proof of Theorem A]
  By Lemma~\ref{LmmLinesInProjectiveSpace}, the variety of lines in $\mathbb P^n$ passing through the given point $y$ can be identified with $H$  and \( \Sigma \subseteq \mathbb{P}^{n-1}\cong H \). Let $L\subset \mathbb P^n$ be a line passing through $y$. By Lemma~\ref{LmmExpressionDefiningSigmaInH}, the line $L$ lies in $Y$ if and only if the equations \( f_i(x_1, \dots, x_n) = 0 \) hold for \( i = m, \dots, d \). This implies that \( \Sigma \) can be viewed as a subvariety of $\mathbb P^{n-1}$, defined by \( d - m + 1 \) homogeneous polynomials \( f_i \) with degree $d_i = i$ for $i = m, m+1, \ldots, d$.

    By Proposition~\ref{PropClassicalResult}, when \( d - m + 1 \leq n - 1 \), namely, when \( d \leq n + m - 2 \), the variety $\Sigma$ is smooth of dimension $ n + m - 2 - d$ of degree $m\cdot (m+1)\ldots d = \frac{d!}{(m-1)!}$. This proves (i). The assertion in (ii) also follows, since when $ d = n + m - 2$, the dimension of $\Sigma$ is $0$. As a set $\Sigma$ is a finite set of $\frac{d!}{(m-1)!}$ points. By Proposition~\ref{PropClassicalResult} again, when $d < n + m - 2$, the variety $\Sigma$ is irreducible. This proves (iii).
\end{proof}
\section{Singularities on Intersections with Linear Spaces}
\label{sec:singularities}
In this section, let us dive into the context of Theorem B. Let $r, n$ be positive integers satisfying the numerical condition
\[
 n + 1 = \binom{ r + 3 }{2}.
\]
Let $Y\subset \mathbb P^n$ be a general smooth cubic hypersurface of dimension $n - 1$ defined by $f=0$. Let $H_{r + 1}$ be a dimension $r+1$ general linear subspace in $\mathbb P^n$ that is tangent to $Y$ along a linear subspace $P$ of dimension $r$ (in particular, $P$ is contained in $Y$). Let $H_{2r + 1}$ be a general linear subspace of dimension $2r+1$ that contains $H_{r + 1}$ as a subspace, and let $Y_{2r} = H_{2r+ 1} \cap Y$. Clearly, the dimension of $Y_{2r}$ is $2r$.
\subsection{Proof of Theorem B (i)}
\begin{lemma}\label{LmmSingY_2rInH_r+1}
    With notations and generality assumptions as above, the singular locus of $Y_{2r}$ is contained in $H_{r+1}$.
\end{lemma}
\begin{proof} Let $|D|$ be the linear system of $Y$ containing the hyperplane sections of $Y$ that contain $H_{r + 1}\cap Y$. Then $H_{2r + 1}$ is the intersection of $n - 2r - 1$ general elements $L_1, \ldots, L_{n - r - 1}$ in $|D|$. It is not hard to see that the base locus of $|D|$ is exactly $H_{r+1}\cap Y$. Let us prove by induction on $i$ that for general elements $L_1, \ldots, L_i$ in $|D|$, the singular locus of the intersection $L_1\cap \ldots \cap L_i$ is contained in $H_{r+1}$. The base case $i=1$ is the simplified version of the Bertini's theorem (\autoref{bertini1}), that is, the singular locus of $L_1$ is contained in the base locus of the linear system $|D|$ which is $H_{r+1}\cap Y$. Now suppose that the singular locus of $L_1\cap\ldots\cap L_i$ is contained in $H_{r+1}\cap Y$. The singular locus of $(L_1\cap \ldots \cap L_i)\cap L_{i+1}$ is contained in the singular locus of $(L_1\cap \ldots \cap L_i)$ and the base locus of $L_{i+1}$, both of which is contained in $H_{r+1}$. Hence, the singular locus of $Y_{2r} = L_1\cap \ldots \cap L_r$ is contained in $H_{2r+1}$, as desired.
\end{proof}

The above lemma suggests us to understand the geometry of $H_{r+1}\cap Y$, where the singular locus of $Y_{2r}$ is located. We have 
\begin{lemma}
    As a set, $H_{r + 1} \cap Y$ is a union of $2$ linear subspaces of dimension $r$, one of which is $P$.
\end{lemma}
\begin{proof} By the Bézout's theorem, the intersection $H_{r+1}\cap Y$ is a cubic hypersurface in $H_{r+1}$. But $$H_{r+1}\cap Y$$ already contains a double linear subspace $P$, since $H_{r+1}$ is tangent to $Y$ along $P$. Since the degree of $H_{r+1}\cap Y$ is $3$, there is a residual linear subspace $P'$ in the intersection $H_{r+1}\cap Y$. Hence, as a set $H_{r+1}\cap Y = P\cup P'$.
\end{proof}

Let $P$, $P'$ be the two linear subspaces of dimension $r$ that appear in $H_{r+1}\cap Y$. A priori, the two linear subspaces $P$ and $P'$ can be identical. However, a general choice of the linear subspace $H_{r+1}$ would ensure $P\neq P'$. Hence, without loss of generality, we may assume that $H_{2r+1}$ is defined by $x_{2r+2} = \ldots = x_n = 0$, that $H_{r+1}$ is defined by $x_{r+1} = \ldots = x_{2r+1}= 0$, that $P$ is defined by $x_{r+1}=x_{r+2}=\ldots=x_{2r+1}=0$, and that $P'$ is defined by $x_0=x_{r+2}=\ldots=x_{2r+1} = 0$. Using the coordinates as above, we have

\begin{lemma}\label{LmmExpressionOfY_2r}
    The defining equation of $Y_{2r}\subset H_{2r+1}\cong\mathbb P^{2r+1}$ is represented as 
    \begin{align*}
       & f(x_0, \ldots, x_{2r + 1})\\
       = &  x_{r + 1}^2 x_0 + x_{r + 2} Q_1(x_0, \ldots, x_{2r + 1}) + \ldots + x_{2r + 1}Q_{r}(x_0, \ldots, x_{2r + 1}), 
    \end{align*}
    where $Q_i(x_0, \ldots, x_{2r + 1})$ are homogeneous quadratic polynomials.
\end{lemma}
\begin{proof}
    Let $\hat{f}$ denote the restriction of $f$ to $H_{2r+1}$, then $\hat{f}$ is the defining polynomial of $Y_{2r}\subset\mathbf{P}^{2r+1}$. By definition, we have \[
    \hat{f}(x_0,\cdots,x_{2r+1})=f(x_0,\cdots,x_{2r+1},0,\cdots,0).
    \]
With some linear re-indexing, the tangency condition forces both $\hat{f}$ and $\frac{\partial \hat{f}}{\partial x_{r+1}}$ to vanish along ${P}\subset H_{r+1}$. Thus, $\hat{f}_|{}_{P}$ must contain the term $x_{r+1}^2$, yielding 
    \[
    \hat{f}_|{}_P=x_{r+1}^2\cdot L(x_0,\cdots\,x_r),
    \]
    where $L$ is a linear function.
    Again, $\hat{f}$ vanishes along $H_{r+1}\supset P^\prime$ defined by $\{x_0=0;x_{r+2}=\cdots=x_{2r+1}\}$. It must admit the locus $\{x_0=0\}$ on $P$, which yields  
    \[
    \hat{f}_|{}_{P\cup P^\prime}=0\implies \hat{f}_|{}_{P}=x_{r+1}^2x_0
    \]
The terms of $\hat{f}(x_0, \ldots, x_{2r + 1})$ readily have the form
    \[
        x_{r + 1}^2 x_0 + x_{r + 2} Q_1(x_0, \ldots, x_{2r + 1}) + \ldots + x_{2r + 1}Q_{r}(x_0, \ldots, x_{2r + 1}),     
    \]   
where $Q_i(x_0, \ldots, x_{2r + 1})$ are as desired.
    
    


\end{proof}

We are finally ready to prove (i) of Theorem B.
\begin{proof}[Proof of Theorem B (i)] Let $x = [ x_0: \ldots: x_{2r+1}]$ be a singular point of $Y_{2r}\subset H_{2r+1} \cong \mathbb P^{2r+1}$. By Lemma~\ref{LmmSingY_2rInH_r+1}, the point $x$ is contained in $H_{r+1}$, i.e., $x_{r+2} = \ldots = x_{2r+1} = 0$. On the other hand, by applying the Jacobian criterion to the expression of $Y_{2r}\subset \mathbb P^{2r+1}$ as described in Lemma~\ref{LmmExpressionOfY_2r}, we find that $x_{r+1} = 0$ and $Q_1(x) = \ldots = Q_r(x) = 0$. We also know that any point $x\in \mathbb P^{2r+1}$ satisfying the conditions $x_{r+1} = x_{r+2} = \ldots = x_{2r+1} = 0$ and $Q_1(x) = \ldots = Q_r(x) = 0$ is a singular point of $Y_{2r}$. The above conditions give a complete intersection of $2r+1$ hypersurfaces in $\mathbb P^{2r+1}$, and the total degree is $2^r$. Hence, as a set, the singular locus of $Y_{2r}$ is a finite set of $2^r$ points.
\end{proof}

\subsection{Proof of Theorem B (ii)}
\begin{lemma}
    The dimension of $\Sigma_y$ is $2r-2$.
\end{lemma}
\begin{proof}
    Similar to Lemma~\ref{LmmExpressionOfY_2r}, taking $m=3,d=0$, we decompose $f$ as
        \begin{align*}
            &f(x_0,x_1,\dots,x_n)\\
            =&x_0^{0}f_3(x_1,\dots,x_n)+
            x_0^{1}f_{2}(x_1,\dots,x_n)+
            x_0^2f_1(x_1,\dots,x_n)+x_0^3.
        \end{align*}
        The last two terms are zero since $y$ is a simple double point. We have $f=f_3+x_0f_2$. A line $l\subset H_{2r+1}\cong\mathbb{P}^{2r+1}$ corresponding to a point $[x_1:...:x_{2r+1}]\in \mathbb{P}^{2r}$, included in $Y_{2r}$, then we have
$f_3(x_1,...,x_{2r+1})=0$ and $f_2(x_1,...,x_{2r+1})=0$. Hence $\Sigma_y$ is a subvariety in $\mathbb{P}^{2r}$ defined by $f_3$ and $f_2$. Hence $\Sigma_y$ is of dim $2r-2$.
\end{proof}

\begin{lemma}
    The singular locus of $\Sigma_y$ is of dimension $\leq 2r - 4$.
\end{lemma}
\begin{proof} First we consider the case $r=2$. We have proved that a point on $\Sigma_y$ is singular only if $x_3=x_4=x_5=0$. Since $\Sigma_y$ is defined by $q=g=0$, a point $(s, t, 0, 0, 0)$ on $\Sigma_y$ is singular iff the rank of 
$$\begin{pmatrix}
\frac{\partial q}{\partial x_1} & \frac{\partial q}{\partial x_2} & ...&\frac{\partial q}{\partial x_5}\\
\frac{\partial g}{\partial x_1} & \frac{\partial g}{\partial x_2} & ...&\frac{\partial g}{\partial x_5}
\end{pmatrix}$$
is less than $2$.

Since $$q=x_3^2+\sum_{i=1}^5(2a_{i0}^1x_ix_4+2a_{i0}^2x_ix_5)$$ and $$g=x_4(\sum_{i, j\geq1}^5a_{ij}^1x_ix_j)+x_5(\sum_{i, j\geq1}^5a_{ij}^2x_ix_j),$$
we deduce that $(s, t, 0, 0, 0)$ on $\Sigma_y$ is singular iff
$$\det\begin{pmatrix}
2a_{12}^1st & 2a_{12}^2st\\
2a_{10}^1s+2a_{20}^1t&2a_{10}^2s+2a_{20}^2t
\end{pmatrix}=0.$$
Solving the equation, we deduce that $\Sigma_y$ has exactly $3$ singular points. In particular, the singular locus of $\Sigma_y$ is of dimension $0\leq0=2r-4$.

When $r\geq3$, the singular points in $\Sigma_r$ correspond to lines in $\mathbb{P}^{r+1}$. Since the modular space of lines through $y$ in $\mathbb{P}^{r+1}$ has dimension $r-1$, the dimension of the singular locus of $\Sigma_r$ is $\leq r-1\leq 2r-4$.  In conclusion, the lemma is true for all $r\geq2$.
\end{proof}

Another ingredient of the proof is the following proposition.
\begin{proposition}\label{PropVariantOfFH}
Let $X\subset \mathbb P^n$ be a complete intersection of irreducible hypersurfaces of dimension $r\geq 2$. Suppose that the singular locus of $X$ is of codimension $\geq 2$ in $X$, then $X$ is irreducible.  
\end{proposition}
\begin{proof}
Let $X = C_1\cup C_2\cup \cdots\cup C_n$ be the irreducible decomposition of $X$. Assuming $X$ is not irreducible, and writing $C_1'=C_2\cup \cdots \cup C_n \neq \emptyset$,  we have $X=C_1\cup C_1'$. Since every regular local ring is factorial, the intersection of $C_1$ and $C_1'$ is contained in the singular locus of $X$, and thus $\dim C_1\cap C_1' \leq r - 2$ by the assumption on the singular locus of $X$. Therefore, we may choose a linear subspace $H\subset \mathbb P^n$ of codimension $r + 1$, such that $H\cap (C_1 \cap C_1') = (H\cap C_1) \cap (H \cap C_1')$ is empty. On the other hand, $\dim (H\cap X) \geq 1$. According to the Fulton-Hansen connectedness theorem (\autoref{connectedness}), the variety $H\cap X$ is connected. However, $H\cap X$ can be decomposed into two disjoint union of nonempty closed subsets $H\cap C_1$ and $H\cap C_1'$. This leads to a contradiction.
\end{proof}
\begin{proof}[Proof of Theorem B (ii)]
    By lemma 4.4, we already have the dimension of $\Sigma_y$ is $2r-2$. By the Proposition 4.6, we set $X=\Sigma_y$. Since $\Sigma_y$ is a complete intersection of irreducible hypersurfaces of $\dim \Sigma_y=2r - 2\geq 2$ (as $r\geq 2$).  By the Lemma 4.5, the singular locus is of codimension $\geq 2$. Therefore, by Proposition~\ref{PropVariantOfFH}, we get that if $r\geq 2$, then $\Sigma_y$ is irreducible of dimension $2r - 2$.
\end{proof}

\end{document}